\documentclass{birkau}

\usepackage{amsmath,amssymb,url}
\usepackage{amssymb}
\usepackage [all]{xy}
\usepackage{mathrsfs}
\usepackage{mathtools}
\usepackage{forest}
\usepackage{stmaryrd}
\usepackage{lmodern}
\usepackage{amsthm}

\numberwithin{equation}{section}

\theoremstyle{plain}
\newtheorem{theorem}{Theorem}[section]
\newtheorem{lemma}[theorem]{Lemma}
\newtheorem{proposition}[theorem]{Proposition}
\newtheorem{corollary}[theorem]{Corollary}
\newtheorem{conjecture}[theorem]{Conjecture}

\theoremstyle{definition}
\newtheorem{definition}[theorem]{Definition}

\newtheorem{remark}[theorem]{Remark}
\newtheorem{example}[theorem]{Example}

\newcommand{\M}{\mathbb{M}}
\newcommand{\Nat}{\mathbb{N}}
\newcommand{\rk}{\mathrm{rk}}

\newcommand{\q}{ \mathfrak{q}}
\newcommand{\Long}{ \mathcal{L}}
\newcommand{\Cat}{\mathbf{C}}

\begin{document}

\title{Growth and density in free magmas\footnote{Under review at \emph{European Journal of Combinatorics}. } }


\author{Carles Card\'o}


\address{Departament de Ci\`encies Bàsiques\\ Universitat Internacional de Catalunya, \\c/ Josep Trueta s/n, Sant Cugat del Vallès 1-3. 08195,  Spain}
\email{ccardo@uic.es}


\maketitle

\begin{abstract}
The density of a subgroupoid with respect to a free groupoid is defined as the asymptotic ratio of their growths. This notion can be interpreted as a generalisation of the index's inverse for groups or as the probability of an element belonging to a subgroupoid. This more combinatorial strategy shows a richer picture of free groupoids than the bare algebraic perspective. We study the growth and density of several subgroupoids of a free groupoid. In addition, some aspects of enumeration related to the Motzkin paths are shown.
\newline
\textbf{keywords: } 
\newline 
\textbf{SMC: } 
\end{abstract}

\section{Introduction} \label{Introduction}  \label{Preliminaries}

Free groupoids are graded by a \emph{length} function. The \emph{growth} of a subgroupoid is the sequence of the ball sizes using the length as a norm. Given a pair of subgroupoids $N$ and $M$, with $N \subseteq M$, of a free groupoid, we define the \emph{density of $N$ respect to $M$} as the asymptotic ratio of the growths of $N$ and $M$. The density, thus, can be interpreted as the likelihood of an element being in $N$ knowing that it is in $M$. 

The equivalent notion of density for other free structures, such as $\Nat$ or $\mathbb{Z}$, do not exhibit more difficulty: the density of $a \mathbb{Z}$ respect to $\mathbb{Z}$ is $1/a$, which coincides with the inverse of the index of the subgroup. Although the index concept can be generalised to semigroups or monoids in several forms \cite{cain2013few}, these pathways are not helpful in more general groupoid theory. Density can be seen as an index of free groupoids. Growth and density show that the structure of free groupoids is more sophisticated than one can suppose at first glance. 

We introduce the main concepts in Section~\ref{Density}, and after some early elementary properties, we visit several examples. The case of finitely generated subgroupoids deserves special attention, Section~\ref{FinitelyGenerated}. \emph{Longitudinal groupoids}, Section~\ref{Longitudinal}, are a sample of subgroupoids for which the density is not defined, although the asymptotic growth ratio can be described.  Last  Section~\ref{CountingSubgroupoid} exhibits a general formula for counting elements of a subgroupoid and shows that \emph{Motzkin numbers} are connected with counting sequences for some subgroupoids. This formula also allows us to propose some subgroupoids with intermediate density.

Following the notation of \cite{cardo2019arithmetic}, $\M_A$ denotes the free groupoid, or the free binar, generated by the set $A$. $\M$ denotes the cyclic free groupoid, and we use the symbol $1$ as a generator and the additive notation. Thus, $\M=\M_{\{1\}}$ and its elements are of the form $1,(1+1), 1+(1+1), (1+1)+1, \ldots$. 
We will use the following abbreviations: $2=1+1$, and
\begin{align*} 
&1_-=1, \,\,\,\,\,\,\,\,\, (n+1)_-=n_-+1,\\
&1_+=1, \,\,\,\,\,\,\,\,\, (n+1)_+=1+n_+,
\end{align*}
for any integer $n\geq 1$.
Given a subset $X\subseteq \M$, $\langle X \rangle$ denotes the least subgroupoid containing $X$. We accept the empty set as groupoid. 
Each subgroupoid $N$ of a free groupoid has a unique minimal generating. The \emph{rank} of $N$ is defined as the cardinality of its minimal generating set.

As described in \cite{cardo2019arithmetic}, $\M$ can be equipped with a product operation defined as the unique operation $\cdot: \M\times \M \longrightarrow \M$ such that for any $x,y,z\in \M$
\begin{align*}
x\cdot 1 &=1\cdot x =x,\\
x\cdot (y+z)&=x\cdot y+x \cdot z.
\end{align*}
$(\M, +, \cdot,1)$ is a \emph{left semi-ringoid} with unity \cite[p.~206]{Rosenfeld1968Algebraic}); that is, the product is associative, whereby $(\M, \cdot, 1)$ is a monoid, and it is left distributive with the sum. 
The product $x\cdot y$ can be regarded as the result of substituting each $1$ by $x$ into $y$. For example $2\cdot 3_+= 2\cdot (1+(1+1))=(2+2(1+1))=(2+(2+2))$. We abbreviate $2^2=2\cdot 2, 2^3= 2\cdot 2^2$, and so forth. 

There is a unique function $\ell: \M_A \longrightarrow \Nat$, called \emph{length}, such that $\ell(a)=1$ for each generator $a \in A$ and $\ell(x+y)=\ell(x)+\ell(y)$ for each $x,y\in \M_A$. In addition, when the free groupoid is cyclic, $\ell$ satisfies that $\ell(xy)=\ell(x)\ell(y)$. 
Given a subset $X$ of $\M_A$, we write $ (X)_n =\{ x\in X \mid \ell(x)=n\}$. For example,
\begin{align*}
(\M)_1&=\{1\},\\
(\M)_2&=\{2\},\\
(\M)_3&=\{3_-,\, 3_+\},\\
(\M)_4&=\{4_-,\, 1{+}3_-,\, 2^2,\, 3_+{+}1,\, 4_+\},\\
(\M)_5&=\{5_-,\, 1{+}4_-,\, 1{+}(1{+}3_-),\, (1{+}3_-){+}1,\, 2{+}3_-, \, 3_-{+}2,\, 2^2{+}1,\\
&\,\,\,\,\,\,\,\,\,\,\,\, 1{+}2^2, \, 2{+}3_+,\, 3_+{+}2,\, 1{+}(3_+{+}1),\,  (3_+{+}1){+}1,\,  4_+{+}1,\,  5_+\}.
\end{align*}

The \emph{counting sequence} of $X$ is defined as $|X|_n=|(X)_n|$. The \emph{growth sequence} of $X$ is defined as $|X|_{\leq n}= \sum_{k=1}^n |X|_n$. 
 It is well-known that $|\M|_n=C_{n-1}$ where $C_n$ are the Catalan numbers defined as 
\begin{align*} C_n=\frac{1}{n+1} { 2n \choose n}  \sim \frac{4^n} {(n+1)\sqrt{\pi n}}\: ,  \end{align*}
for every $n\geq 0$, \cite{Stanley2015}, \cite{crepinvsek2009efficient}, \cite{Dutton1986Comput}, and $\sim$ stands for asymptotic equivalence.  
We will write $c_0=0$ and $c_n=C_{n-1}$ for each $n>0$. Thus, $|\M|_n=c_n$.  
Catalan numbers can also be calculated by the recursive formulas
\begin{align*}
C_{n+1}=\sum _{i=0}^{n}C_{i}\,C_{n-i}\quad \mbox{ or } \quad C_{n+1}={\frac {2(2n+1)}{n+2}}C_{n}\: ,
\end{align*}
for $n\geq 0$ and $C_0=1$. The following inequalities hold \cite{Dutton1986Comput},
\begin{align*} \frac{4^n}{(n+1)\sqrt{\pi n \left(\frac{4n}{4n-1} \right)}}<C_n<\frac{4^n}{(n+1)\sqrt{\pi n \left( \frac{4n+1}{4n} \right) }}\: .
\end{align*}
Since for $n\geq 4$, $\sqrt{\pi n\frac{4n}{4n-1}}\leq n \leq n+1$, these bounds can be  weakened:
\begin{align*} \frac{4^n}{(n+1)^2}<C_n<4^n .
\end{align*}

\section {Density of subgroupoids of a free groupoid} \label{Density}

\begin{definition} Given two subgroupoids, or more in general, two subsets of a free groupoid, $N \subseteq M\subseteq \M_A$, the \emph{density of $N$ respect to $M$} is defined as the limit
\begin{align*}
\delta(N: M)=\lim_{n\rightarrow \infty} \frac{|N|_{\leq n}}{|M|_{\leq n}}\: .
\end{align*}
In what follows, when we write the expression $\delta(N:M)$, we will assume that the above limit exists. 
\end{definition}

\begin{remark} The \emph{Stolz–Cesàro theorem}, see for example \cite[p.~85]{muresan2009concrete}, states that given two sequences of real numbers $A_n$ and $B_n$ and assuming that $B_n$ is strictly monotone and divergent, if the limit $\lim_{n \to \infty } \frac{A_n}{B_n}=L$ exists, then
\begin{align*}
\lim_{n \to \infty} \frac{A_1+\cdots +A_n}{B_1+\cdots +B_n}=L\: .
\end{align*}
In our case, that means that
\begin{align*}
\delta(N: M)=\lim_{n\rightarrow \infty} \frac{|N|_{\leq n}}{|M|_{\leq n}}=\lim_{n\rightarrow \infty} \frac{|N|_{ n}}{|M|_{ n}}\: ,
\end{align*}
provided that the last limit exists. We will use that result in the calculations when necessary without commenting on it. 
\end{remark}

\begin{proposition} For any subgroupoids $N \subseteq N' \subseteq N'' \subseteq \M$,
\begin{itemize}
\item[(i)] $0\leq \delta(N:N')\leq1$;
\item[(ii)] $\delta(N:N)=1$;
\item[(iii)] $\delta(N:N')\cdot \delta(N':N'')=\delta(N:N'')$;
\item[(iv)] $\delta(N:N'')\leq \delta(N':N'')$.

\end{itemize}
\end{proposition}
\begin{proof} All the properties are straightforward from the definition. 
\end{proof}

In view of property (iii), we are interested in densities of the form $\delta(N: \M)$, since the density $\delta(N:M)$ can be calculated from the quotient $\delta(N: \M)/\delta(M: \M)$, provided that the denominator is not zero. 
Let us see some non-trivial properties. We recall, \cite{cardo2019arithmetic}, that given a subgroupoid $N \subseteq \M$, we have the inclusion $Na \subseteq N$, although in general, $Na$ is not a subgroupoid, and that $aN$ is a subgroupoid but in general $aN\not \subseteq N$.

\begin{proposition} Let $N \subseteq M \subseteq \M$ be a pair of subgroupoids and let $a\in \M$,
\begin{itemize}
\item[(i)] $\delta(aN:aM)=\delta(N:M)$;
\item[(ii)] $\delta(Na: Ma)=\delta(N:M)$; 
\item[(iii)] $\delta(Na:N)=\delta(Ma:M)$, provided $\delta(N:M)\not=0$.
\end{itemize}
\end{proposition}
\begin{proof} (i) Firstly we notice that for any subgroupoid $H \subseteq \M$ we have that  $|aH|_{\ell(a)n}=|H|_n$. Then, 
\begin{align*}
\delta(aN: aM)=\lim_{n \to \infty} \frac{|aN|_{\leq n}}{|aM|_{\leq n}}=\lim_{n \to \infty} \frac{|aN|_{\leq \ell(a)n}}{|aM|_{\leq \ell(a)n}}=\lim_{n \to \infty} \frac{|N|_{\leq n}}{|M|_{\leq n}}=\delta(N:M)\: .
\end{align*}
In the second equality, we have used that, since the limit of $|aN|_{\leq n}/|aM|_{\leq n}$ exists, it can be accelerated by the subsequence $|aN|_{\leq \ell(a)n}/|aM|_{\leq \ell(a)n}$. The second statement (ii) is similar to (i). For (iii) consider the inclusions,
\begin{center}
\begin{tabular}{c c c}
  $Na$ & $ \subseteq $ & $N$\\
  \rotatebox[origin=c]{270}{$\subseteq$} &&   \rotatebox[origin=c]{270}{$\subseteq$} \\
  $Ma$ & $ \subseteq $ & $M$
\end{tabular}
\end{center}
That yields two equalities:
\begin{align*}
\delta(Na:N)\delta(N:M)&=\delta(Na:M)\: ,\\
\delta(Na:Ma)\delta(Ma:M)&=\delta(Na:M)\: .
\end{align*}
Equalising the left-hand sides, we obtain
\begin{align*}
\delta(Na:N)\delta(N:M)=\delta(Na:Ma)\delta(Ma:M)\: .
\end{align*}
By (ii), we can cancel the terms $\delta(N:M)$ and $\delta(Na:Ma)$ provided they are not null. 
 \end{proof}

The following are a few early, accessible examples. 

\begin{example} Assuming that $A$ is finite, $|\M_A|_n=c_n |A|^n$. If $A$ is a proper subset of a finite set $B$, then $|A|<|B|$, 
\begin{align*}
\delta(\M_A:\M_B) &=\lim_{n \to \infty} \frac{ |\M_A|_{\leq n}}{|\M_B|_{\leq n}}=\lim_{n \to \infty} \frac{ \sum_{k=1}^n |\M_A|_{k}}{\sum_{k=1}^n |\M_B|_{k}}=\lim_{n \to \infty} \frac{ |\M_A|_{n}}{|\M_B|_{n}}\\
&=\lim_{n \to \infty} \frac{c_n|A|^n}{c_n |B|^n}=\lim_{n \to \infty} \left( \frac{|A|}{|B|} \right)^n=0\: .
\end{align*}
\end{example}

\begin{example} \label{Exemple1} Consider the density of subgroupoids with only one generator, $\langle a \rangle =a \M$. Let us show that $\delta(a\M : \M)=0$, provided $a\not=1$, otherwise $\delta(a\M: \M)=\delta(\M:\M)=1$. We have that
\begin{align*}
|a\M|_n=\begin{cases} c_{\frac{n}{\ell(a)}}  & \mbox{ if } \ell(a)|n, \\
0  & \mbox{ otherwise. } \end{cases}
\end{align*}
Therefore, $ |a\M|_{\ell(a)n} = c_n$. Firstly, we calculate the limit of the subsequence
\begin{align*}
\lim_{n\rightarrow \infty} \frac{|a\M|_{ \ell(a)n}}{|\M|_{\ell(a)n} }=\lim_{n\rightarrow \infty} \frac{c_n}{c_{\ell(a)n}}\: .
\end{align*}
Using the inequalities from Section~\ref{Preliminaries}, $\frac{4^{n-1}}{n^2}<c_n<4^{n-1}$, and since $\ell(a)>1$,
\begin{align*}
\lim_{n\rightarrow \infty} \frac{c_n}{c_{\ell(a)n}}< \lim_{n\rightarrow \infty} \frac{4^{n-1}}{4^{\ell(a)n-1}}n^2=\lim_{n\rightarrow \infty} \left( \frac{1}{4^{\ell(a)-1}}\right)^n n^2=0\: .
\end{align*}
Since the subsequence $|a\M|_{ \ell(a)n}/|\M|_{\ell(a)n}$ tends to zero, and the complementary terms are null, $\lim_{n\rightarrow \infty}  |a\M|_{n}/|\M|_{n}=0$. Therefore, $\delta(a\M: \M)=0$.
\end{example}

\begin{example} A subgroupoid $N \subseteq \M$ it is said to be \emph{finitely maximal} if the set of subgroupoids $M$ such that $N \subseteq M \subseteq \M $ is finite. Some examples are
\begin{align*}
\M, \quad \M \setminus \{ 1\}, \quad  \M \setminus \{1, 2, 2+2\}, \quad\M \setminus \{1, 3_-, 5_+\}\: .
\end{align*}
We studied such subgroupoids in \cite[Theorem~6.10]{cardo2019arithmetic}, proving that $N$ is finitely maximal if and only if its complement $\M \setminus N$ is finite. 
Trivially, finitely maximal subgroupoids of $\M$ have full density. Consider the sets
 \begin{align*}
 \q_n=\{1, 2, 3_+, 4_+, 5_+, \ldots, n_+ \} \mbox{ and } \q=\{ n_+\mid n \in \Nat\}\: .
 \end{align*}
 Both $\M\setminus \q_n$ and $\M \setminus \q$ are subgroupoids with density one. 
However, the subgroupoids $\M\setminus \q_n$ are finitely maximal, while $\M\setminus \q$ is not, which shows that non-finitely maximal subgroupoids with full density exist. 
\end{example}
 
\section{Density of finiteley generated subgroupoids} \label{FinitelyGenerated}

Consider the density of finitely generated subgroupoids of the free groupoid. As we saw in Example~\ref{Exemple1}, proper principal ideals have null density. The following is a more general result. 
\begin{theorem} \label{TheoFG} Let $G$ the minimal generating set of a non-empty subgroupoid $N \subseteq \M$. Set $\lambda(N)=\min \{ \ell(x) \mid x \in G\}$.  If 
$\rk(N)< 4^{\lambda(N)-1}$, then $N$ has null density. 
\end{theorem}\label{FGgroupoids}
\begin{proof} First of all we recall the \emph{multinomial theorem} \cite{knuth1989concrete}, which states that given $m$ real numbers $x_1$, $x_2$, $\ldots, x_m$,
\begin{align*}
\sum _{k_{1}+k_{2}+\cdots +k_{m}=n}{n \choose k_{1},k_{2},\ldots ,k_{m}}x_{1}^{k_{1}}x_{2}^{k_{2}}\cdots x_{m}^{k_{m}}=(x_{1}+x_{2}+\cdots +x_{m})^{n}.
\end{align*}
where the \emph{multinomial coeficients} are defined as
\begin{align*} {n \choose k_{1},k_{2},\ldots ,k_{m}}=\frac {n!}{k_{1}!\,k_{2}!\cdots k_{m}!}\: .
\end{align*}
The substitution of $x_i=1$ for all $i$ into the above expression gives immediately that
\begin{align*}
\sum _{k_{1}+k_{2}+\cdots +k_{m}=n}{n \choose k_{1},k_{2},\ldots ,k_{m}}=m^{n}.
\end{align*}

Let $N$ be the finitely generated subgroupoid $N=\langle a_1, \ldots, a_r \rangle$, with rank $\rk(N)=r\geq 1$ and assume that the set of generators is minimal. We set $\alpha_1=\ell(a_1), \ldots, \alpha_r=\ell(a_r)$, and $\lambda=\lambda(N)=\min \{\alpha_1, \ldots, \alpha_r\}$. 

Think of the elements in $\M$ as planar binary trees. An element in $x \in N$ with length $\ell(x)=n$ can be built by taking some tree $y\in \M$ and hanging a subtree $a_i$ from each leaf of $y$, with $i=1, \ldots, r$. 
Therefore, $n=q_1\alpha_1+ \cdots + q_r \alpha_r$, where each $q_i$ is the number of occurrences of the subtree $a_i$ in $x$. Some of these $q_i$ can be $0$, but at least there is some $q_i>0$. Notice that $\ell(y)=q_1+\cdots +q_r$ and that there are 
\begin{align*}
{ q_1+\cdots +q_r \choose q_1, \ldots, q_r}
\end{align*}
ways to hang the generators in $y$. Since there are $c_{\ell(y)}=c_{q_1+\cdots+q_r}$ possible trees $y$ of length $\ell(y)$, we can calculate $|N|_n$ by the expression
\begin{align*}
|N|_n=\sum_{\mathcal{A}} { q_1+\cdots +q_r \choose q_1, \ldots, q_r} c_{q_1+\cdots+q_r}\: ,
\end{align*}
where the sum runs over $(q_1, \ldots, q_r) \in \mathcal{A}$ defined as
\begin{align*}
\mathcal{A}=\{ (q_1, \ldots, q_r) \in (\Nat \cup \{0\})^r \,\,\mid\,\, \alpha_1 q_1+\cdots + \alpha_r q_r=n, \,\,\,\, q_1+\cdots + q_r>0 \}\: .
\end{align*}
Now we observe that
\begin{align*}
n=\alpha_1q_1 + \cdots + \alpha_r q_r\geq \lambda q_1+ \cdots + \lambda q_s=\lambda (q_1 +\cdots + q_r)\: ,
\end{align*}
therefore,
\begin{align*}
q_1 +\cdots + q_r \leq \lfloor \frac{n}{\lambda} \rfloor\: ,
\end{align*}
where $\lfloor x \rfloor$ is the floor function. That means that if the set of integers $q_1, \ldots, q_r$ is a solution of the equation $\alpha_1 q_1+\cdots + \alpha_r q_r=n$, then it is also a solution of the inequation $q_1+\cdots +q_r \leq \lfloor \frac{n}{\lambda} \rfloor $. This serves to bound the number $|N|_n$ as
\begin{align*}
|N|_n&=\sum_{\mathcal{A}} { q_1+\cdots +q_r \choose q_1, \ldots, q_r} c_{q_1+\cdots+q_r}  \leq \sum_{q_1+\cdots + q_s \leq \lfloor \frac{n}{\lambda} \rfloor } { q_1+\cdots +q_r \choose q_1, \ldots, q_r} c_{q_1+\cdots+q_r} \\
&= \sum_{k=1}^{\lfloor \frac{n}{\lambda} \rfloor } \,\,\, \sum_{q_1+\cdots + q_r =k} { k \choose q_1, \ldots, q_r} c_{k} = \sum_{k=1}^{\lfloor \frac{n}{\lambda} \rfloor } c_k \sum_{q_1+\cdots + q_r =k} { k \choose q_1, \ldots, q_r} \\
&= \sum_{k=1}^{\lfloor \frac{n}{\lambda} \rfloor } c_k r^k \leq \sum_{k=1}^{\lfloor \frac{n}{\lambda} \rfloor } 4^{k-1} r^k  = \frac{1}{4} \sum_{k=1}^{\lfloor \frac{n}{\lambda} \rfloor } (4r)^k = \frac{1}{4} \cdot \frac{(4r)^{\lfloor \frac{n}{\lambda} \rfloor +1 }-4r}{ 4r-1}\\
& = \frac{r}{4r -1} \cdot \left( (4r)^{\lfloor \frac{n}{\lambda} \rfloor} -1 \right) \leq \frac{r}{4r -1} \cdot \left( (4r)^{ \frac{n}{\lambda} } -1 \right),
\end{align*}
where we have used the multinomial theorem and the inequality $c_k\leq 4^{k-1}$ in the third line, and that $\lfloor \frac{n}{\lambda} \rfloor  \leq \frac{n}{\lambda}$ in the last inequality. With all this, we have
\begin{align*}
\frac{|N|_n}{|\M|_n} \leq \frac{r}{4r -1} \cdot \frac{ (4r)^{\frac{n}{\lambda} } -1 }{c_n}\,.
\end{align*}
Now we use the inequality, Section~\ref{Preliminaries}, for $n\geq 4$:
\begin{align*}
\frac{4^{n-1}}{n^2}<c_n\: ,
\end{align*}
and then:
\begin{align*}
\frac{|N|_n}{|\M|_n} \leq \left( \frac{(4r)^{ \frac{n}{\lambda} } -1 }{4^n} \right) \left( \frac{r}{4r -1}  4n^2 \right).
\end{align*}
Recall that by hypothesis $r<4^{\lambda -1}$, which means that
\begin{align*}
\frac{(4r)^{\frac{1}{\lambda}}}{4}<1\: ,
\end{align*}
and then, the limit
\begin{align*}
\lim_{n \to \infty } \frac{(4r)^{ \frac{n}{\lambda} } -1 }{4^n} =0\: ,
\end{align*}
dominates the above general limit. That is, $\lim_{n \to \infty } |N|_n/|\M|_n=0$, which means that $\delta(N: \M)=0$. 
 \end{proof}

\begin{corollary} Let $N$ be a proper subgroupoid of $\M$. 
\begin{itemize}
\item[(i)] If $\rk(N)\leq 3$, then it has null density.
\item[(ii)] If $2$ is not a generator of $N$ and $\rk(N)\leq 15$, then it has null density.
\item[(ii)] If $2, 3_-, 3_+$ are not generators of $N$ and $\rk(N)\leq 63$, then it has null density.
\end{itemize}
\end{corollary}
\begin{proof} If $N$ is empty, then the statement is trivial in all the cases. (i) If $N$ is proper, then $N\not= \M$, whereby $1$ is not a generator of $N$ and  hence, $\lambda (N)\geq 2$. By Theorem~\ref{FGgroupoids},  if $\rk(N)<4^{\lambda(N)-1}\leq 4^{2-1}=4$, then $N$ has null density. (ii) is similar to (i) with, in addition, $2$ is not a generator, whereby $\lambda (N)\geq 3$. For (iii) we use $\lambda(N)\geq 4$. 
 \end{proof}

\begin{example} We can also obtain the density of a subgroupoid, although the hypothesis of Theorem~\ref{FGgroupoids} does not hold. Consider the finitely generated subgroupoid
\begin{align*}
N=\langle 3_+, 2\cdot 2, 2\cdot 3_+, 2\cdot 4_+, \ldots, 2\cdot 16_+\rangle\: .
\end{align*}
We can see that $\rk(N)=16$ and $\lambda(N)=3$, and then,
\begin{align*}
16=\rk(N)\not< 4^{\lambda(N)-1}=4^2=16\: .
\end{align*}
However, $N \subseteq \langle 3_+, 2\rangle$, and therefore $\delta(N, \M) \leq \delta( \langle 3_+, 2 \rangle: \M)=0$.
\end{example}

\begin{proposition} \label{PropaM} Let $N$ be any subgroupoid of $\M$ and let $a \in \M$. If $a\not=1$, then $\delta(aN: \M)=0$.
\end{proposition}
\begin{proof} Since $N \subseteq \M$, $aN \subseteq a\M$ and $\delta(aN: \M)=\delta (aN: a\M)\cdot \delta( a\M: \M)= \delta (aN: a\M)\cdot 0=0$. 
 \end{proof}

\begin{remark} \label{RemarkConjecture} Proposition~\ref{PropaM} implies that there are subgroupoids non-finitely generated with null density. Let $N$ be a non-finitely generated subgroupoid, say for example $N =\langle 2, 3_+, 4_+, 5_+, \ldots \rangle$, and let $a \in \M$, $a\not=1$. $aN$ is still a non-finitely generated subgroupoid,
\begin{align*}
aN= \langle a2, a3_+, a4_+, a5_+, \ldots \rangle\: ,
\end{align*}
however, it has a null density. 
\end{remark}

\begin{conjecture} \label{Conjecture} All the finitely generated proper subgroupoids of the cyclic free groupoid have null density.
\end{conjecture}
In the last Section~\ref{CountingSubgroupoid}, we will do some empirical calculations that reinforce that conjecture. 
If some counterexample exists, this should have more than three generators and, by Proposition~\ref{PropaM}, they should be coprime in the product operation in $\M$.

\section{Growth of longitudinal subgroupoids} \label{Longitudinal}
An element of a free groupoid is often represented as a binary tree.  
What is the probability that a random binary tree has an even number of leaves? Let us show that the question has no answer, at least in the sense of the probability understood as density. 
A \emph{longitudinal subgroupoid} $L \subseteq \M$ is a subgroupoid of the form  $L=\ell^{-1}(N)$, where $N$ is a subsemigroup of $\Nat$. Longitudinal subgroupoids are not finitely generated, and they form a complete sublattice of the lattice of subgroupoids of the cyclic free groupoid \cite{cardo2019arithmetic}. We will see that the density is not defined for such groupoids, but we can describe its asymptotic behaviour well. Consider firstly the following specific kinds of longitudinal subgroupoids,
\begin{align*}
\Long_p=\{ x \in \M \mid p|\ell(x)\}=\ell^{-1}(p \Nat)\: .
\end{align*}

\begin{lemma} \label{Lemapreviprevi} For any integer $p>0$,
\begin{align*}
\lim_{n \to \infty} \frac{c_{np}}{c_p+c_{2p}+\cdots +c_{np}}=1-\frac{1}{4^p}\: .
\end{align*}
\end{lemma}
\begin{proof} On the one hand,
\begin{align*}
\lim_{n \to \infty} \frac{\sum_{k=1}^n c_{np}-c_{(n-1)p} }{\sum_{k=1}^n c_{np}}= \lim_{n\to \infty} \frac{c_{np}-c_{(n-1)p}}{c_{np}}=1- \lim_{n\to \infty} \frac{c_{(n-1)p}}{c_{np}}=1-\frac{1}{4^p}\: ,
\end{align*}
where we have used the asymptotic approximation $c_n \sim \frac{4^{n-1}} {n\sqrt{\pi (n-1)}}$. On the other hand, since $c_0=0$,
\begin{align*}
\lim_{n \to \infty} \frac{\sum_{k=1}^n c_{np}-c_{(n-1)p} }{\sum_{k=1}^n c_{np}}&=\lim_{n \to \infty} \frac{c_{np}- c_0}{c_p+c_{2p}+\cdots +c_{np}}\\
&=\lim_{n \to \infty} \frac{c_{np}}{c_p+c_{2p}+\cdots +c_{np}}=1-\frac{1}{4^p}\: .
\end{align*}
\end{proof}

\begin{lemma} \label{LemmaLongitudinal} For any $p\geq 1$,
\begin{align*}
\frac{|\Long_p|_{\leq n}}{|\M|_{\leq n}} \sim \frac{3}{\left(1-\frac{1}{4^p} \right)4^{1+(n \,\,\mathrm{mod}\, p)}}\: .
\end{align*}

\end{lemma}
\begin{proof} We have to prove that if $n=pk+r$, with $r=0, 1, \ldots, p-1$, then:
\begin{align*}
\lim_{k\to \infty} \frac{|\Long_p|_{\leq pk+r}}{|\M|_{\leq pk+r}} = \frac{3}{4^{r+1}\left(1-\frac{1}{4^p} \right)}\: .
\end{align*}
First, we prove the case for $r=0$. We notice that
\begin{align*}
|\Long_p|_n=\begin{cases} c_n & \mbox{ if } p | n, \\
0 & \mbox{ otherwise,} \end{cases}
\end{align*}
We sum in groups of $p$ elements the terms of the successions $|\Long_p|_n$ and $|\M|_n$ as follows:
\begin{align*}
|\Long_p|_n &=\overbrace{0,0,\ldots,0\,,c_p}^{A_1}, \overbrace{\,0,\,\,0,\,\, \ldots, \,\,0,\,\, c_{2p}\,}^{A_2}, \,\, \overbrace{\,0, \,\,0,\,\, \ldots, \,\,0, \,\, c_{3p}}^{A_3},\,\,\ldots\\
|\M|_n &=\underbrace{c_1,c_2,\ldots \,,c_p}_{B_1},\, \underbrace{c_{p+1},c_{p+2}, \ldots, c_{2p}}_{B_2}, \underbrace{c_{2p+1},c_{2p+2}, \ldots, c_{3p}}_{B_3},\ldots
\end{align*}
That is,
\begin{align*}
A_k & =c_{pk}\: ,\\
B_k & =c_{pk-(p-1)}+c_{pk-(p-2)}+\cdots +c_{pk-2}+c_{pk-1}+c_{pk}\: .
\end{align*}
If the limit $\lim_{k \to \infty} A_k/B_k$ exists, then, by the Stolz–Cesàro theorem,
\begin{align*}
\lim_{k\to \infty} \frac{|\Long_p|_{\leq pk}}{|\M|_{\leq pk}}=\lim_{k\to \infty} \frac{\sum_{s=1}^k A_s}{\sum_{s=1}^k B_s}=\lim_{k\to \infty} \frac{A_k}{B_k}\: .
\end{align*}
Therefore, we only need to calculate the last limit. Using the recurrence of the Catalan numbers $c_{n-1}=\frac{n+1}{4n-2}c_n$, we can rewrite $B_k$ as
\begin{align*}
B_k=c_{pk} \Big( D_{pk-p+1}\cdots D_{pk-1}D_{pk} +  \cdots + D_{pk-1}D_{pk}  +  D_{pk}  +1 \Big)\: ,
\end{align*}
where $D_n=\frac{n+1}{4n-2}$. Now we take the quotient
\begin{align*}
\frac{A_k}{B_k}=\frac{1}{D_{pk-p+1}\cdots D_{pk-1}D_{pk} +  \cdots + D_{pk-1}D_{pk}  + D_{pk}  +1 }\: .
\end{align*}
Since $\lim_{k \to \infty } D_{pk-s}= \frac{1}{4}$ for each $p$ and each $s=0,1 \ldots, p-1$, when we take the limit, we get
\begin{align*}
\lim_{k \to \infty } \frac{A_k}{B_k}=\frac{1}{1+\frac{1}{4}+\frac{1}{4^2}+\cdots +\frac{1}{4^{p-1}}}=\left( \frac{1-\frac{1}{4^p}}{1-\frac{1}{4}}\right)^{-1}=\frac{3}{4 \left( 1- \frac{1}{4^p} \right)}\: .
\end{align*}
Now we see the general case:
\begin{align*}
\lim_{k\to \infty} \frac{|\Long_p|_{\leq pk+r}}{|\M|_{\leq pk+r}} &=\left( \lim_{k\to \infty} \frac{|\M|_{\leq pk+r}}{|\Long_p|_{\leq pk+r}}\right)^{-1}\\
&=\left( \lim_{k\to \infty} \frac{|\M|_{\leq pk} + (c_{pk+1}+c_{pk+2}+\cdots+c_{pk+r})}{|\Long_p|_{\leq pk+r}}\right)^{-1}\\
&=\left( \lim_{k\to \infty} \frac{|\M|_{\leq pk} + (c_{pk+1}+c_{pk+2}+\cdots+c_{pk+r})}{|\Long_p|_{\leq pk}}\right)^{-1}\\
&=\left( \lim_{k\to \infty} \frac{|\M|_{\leq pk} }{|\Long_p|_{\leq pk}} + \lim_{k\to \infty} \frac{ c_{pk+1}+c_{pk+2}+\cdots+c_{pk+r}}{|\Long_p|_{\leq pk}}\right)^{-1}\\
&=\left( \left(\frac{3}{4 \left( 1- \frac{1}{4^p} \right)}\right)^{-1} + \lim_{k\to \infty} \frac{ c_{pk+1}+c_{pk+2}+\cdots+c_{pk+r}}{c_p+c_{p2}+\cdots +c_{pk}}\right)^{-1}.
\end{align*}
Calculate the last limit. Using the recurrence $c_{n+1}=\frac{4n-2}{n+1}c_n$ we can rewrite the numerator as
\begin{align*}
&\, c_{pk+1}+c_{pk+2}+\cdots+c_{pk+r}\\
=&\, c_{pk}E_{pk}+ c_{pk}E_{pk}E_{pk+1}+\cdots +c_{pk}E_{pk}E_{pk+1}\cdots E_{pk+r-1}\\
=&\, c_{pk} \left( E_{pk}+ E_{pk}E_{pk+1}+\cdots +E_{pk}E_{pk+1}\cdots E_{pk+r-1}\right),
\end{align*}
where $E_n=\frac{4n-2}{n+1}$. Therefore, the limit is
\begin{align*}
&\lim_{k\to \infty} \frac{ c_{pk+1}+c_{pk+2}+\cdots+c_{pk+r}}{c_p+c_{p2}+\cdots +c_{pk}}\\
=& \left( \lim_{k\to \infty}  \frac{c_{pk}}{c_p+c_{p2}+\cdots+c_{pk}} \right) \cdot \left( \lim_{k\to \infty} E_{pk}+ E_{pk+1}+\cdots +E_{pk}E_{pk+1}\cdots E_{pk+r-1} \right) \\
=&\left( 1- \frac{1}{4^p}\right) \cdot (4+4^2+\cdots 4^r)=\left( 1- \frac{1}{4^p}\right) \cdot  \frac{4}{3}(4^r-1)\: ,
\end{align*}
where we have used the limit of Lema~\ref{Lemapreviprevi} and that $\lim_{k \to \infty } E_{pk+s}=4$ for any $p $ and $s=0, 1, \ldots, r-1$. Now we can resume the main limit:
\begin{align*}
\lim_{k\to \infty} \frac{|\Long_p|_{\leq pk+r}}{|\M|_{\leq pk+r}} &=\left( \left(\frac{3}{4 \left( 1- \frac{1}{4^p} \right)}\right)^{-1} +  \left( 1- \frac{1}{4^p}\right) \cdot  \frac{4}{3}(4^r-1)       \right)^{-1} \\
&=\frac{3}{4^{r+1}\left(1-\frac{1}{4^p} \right)}\: . \qedhere
\end{align*}
\end{proof}

\begin{example} Let us show some particular cases. For $p=2$,
\begin{align*}
\frac{|\Long_2|_{\leq n}}{|\M|_{\leq n}} \sim \begin{cases} \frac{4}{5} & \mbox{ if }n \mbox{ is even,} \\ \frac{1}{5} & \mbox{ if }n \mbox{ is odd.}  \end{cases}
\end{align*}
The ratio oscillates between two values as $n$ goes to infinity, and it does not make meaning to interpret it as the probability of an element being in $\Long_2$, or what is the same, as the probability of a random binary tree having an even number of leaves. 
The same occurs for $p=3$,
\begin{align*}
\frac{|\Long_3|_{\leq n}}{|\M|_{\leq n}} \sim \begin{cases} \frac{16}{21} & \mbox{ if } n\equiv 0 \mod 3, \\ \frac{4}{21} & \mbox{ if } n\equiv 1 \mod 3,  \\ \frac{1}{21} & \mbox{ if } n\equiv 2 \mod 3.  \end{cases}
\end{align*}
However, notice that the mean of the set of the asymptotic values for a fixed $p$ is
\begin{align*}
\frac{1}{p} \sum_{r=0}^{p-1} \lim_{k \to \infty} \frac{|\Long_p|_{\leq pk+r}}{|\M|_{\leq pk+r}}=\frac{1}{p}\sum_{r=0}^{p-1} \frac{3}{4^{r+1}\left(1-\frac{1}{4^p} \right)}=\frac{1}{p}\: .
\end{align*}
\end{example}

Now, let us prove that the growth of any longitudinal subgroupoid is asymptotically the same as the growth of $\Long_p$ for some positive integer $p$. Given a subsemigroup of $\Nat$, say $\langle A\rangle$, where $A$ is any subset of $\Nat$, we can write
\begin{align*}
\langle A \rangle = \gcd(A) \cdot \left\langle \frac{A}{\gcd(A)} \right\rangle.
\end{align*}
We set $A'=\frac{A}{\gcd(A)}$. The set $\langle A' \rangle$ pertains to a class of subsemigroups called \emph{numerical semigroup}. See \cite{Finch2009Monoids}, \cite{delgado2013numerical}, \cite{rosales2004numerical} for a survey and some aspects of the issue. It is proved that there is a number, called the \emph{Frobenius number} and denoted by $F(A')$, such that $n \in \langle A' \rangle$ for any $n > F(A')$. In particular, when a numerical semigroup has rank two,
\begin{align*}
F(\{ a, b \})=ab-a-b\: ,
\end{align*}
meanwhile, for greater ranks, there is no closed formula.
\begin{theorem} \label{Ultima} If $L$ is a non-empty longitudinal subgroupoid, then there is an integer $p>0$ such that
\begin{align*}
\frac{|L|_{\leq n}}{|\M|_{\leq n}} \sim \frac{3}{\left(1-\frac{1}{4^p} \right)4^{1+(n \,\,\mathrm{mod}\, p)}}\: .
\end{align*}
\end{theorem}
\begin{proof}  We set $L=\ell^{-1}(N)$, with $N$ a non-empty subsemigroup of $\Nat$, $N=\langle A \rangle$, $p=\gcd(A)$, and $A'=\frac{A}{p}$. The previous comments say that each subsemigroup $N \subseteq \Nat$ behaves as $p \Nat$ for $n \in N$ sufficiently large, or more exactly, for any $n>p F(A')$. This means that $(\ell^{-1}(N))_n=(\ell(p\Nat))_n=(\Long_p)_n$ for any $n>p F(A')$. Thus, we can use the same calculations of the proof of Lemma~\ref{LemmaLongitudinal}.
\end{proof}

\section{Counting elements of a subgroupoid} \label{CountingSubgroupoid}

In the proof of Theorem~\ref{TheoFG} we saw how to count the number of elements of a finitely generated subgroupoid by a multinomial sum.
Now, we show a recursive formula that permits the calculation of non-finitely generated subgroupoids. That will allow us to propose some subgroupoids with a possible intermediate density.

\begin{definition}
Given a sequence of integers $a_n$ we define the \emph{Catalan transform} $\Cat(a_n)=b_n$ as the sequence defined by $b_1=a_1$ and
\begin{align*}
b_n= a_n + \sum_{\underset{0<i,j<n}{i+j=n}} b_ib_j\: .
\end{align*}
\end{definition}
The first elements of the $\Cat$ transform of $a_n$ are
\begin{align*}
b_1&=a_1\: ,\\
b_2&=a_1^2+a_2\: ,\\
b_3&=2a_1^3+2a_1a_2+a_3\: ,\\
b_4&=4a_1^4+a_1^2+4a_1^2a_2+2a_1a_3+a_2+a_4\: ,\\
b_5&=12a_1^5+2a_1^3+16a_1^3a_2+4a_1^2a_3+4a_1a_2^2+2a_2a_3+a_5\: ,\\
\vdots&
\end{align*}
A pair of examples of the $\Cat$ transform:
$$\Cat(1,0,0,0,0,\ldots)=c_n\: , \quad \Cat(-1,0,0,0,0,\ldots)=(-1)^nc_n\: .$$
The last is a particular case of the following nice property,
$$\Cat(\alpha^n \cdot a_n)=\alpha^n \cdot \Cat(a_n)\: ,$$
for any real number $\alpha$.
If $\Phi(x)=\sum_{n=0}^\infty a_nx^n$ is the ordinary generating function of the sequence $a_n$, where we assume $a_0=0$, and if $\Psi(x)$ is the generating function of $\Cat(a_n)$, then, since the generating function of a convolution product turns into a product of generating functions, we have
\begin{align*}
\Psi(x)=\Psi(x)^2+\Phi(x)\: , \quad \mbox{ that is, } \quad \Psi(x)=\frac{1-\sqrt{1-4\Phi(x)}}{2}\: ,
\end{align*}
where the minus sign for the root is conveniently chosen for the same reason as in the case of the generating function of Catalan numbers \cite[p.~4]{Stanley2015}.
\begin{theorem} If $G$ the minimal generating set of the subgroupoid $N\subseteq \M$, then $|N|_n = \Cat(|G|_n)$.
\end{theorem}
\begin{proof} We need an easy result from \cite[Lemma~4.6]{cardo2019arithmetic} 
which states that if $x+y \in N$, then either $x,y \in N$ or $x+y \in G$. 
Hence, $(N)_n$ is formed by sums  $x+y$ such that $x,y \in N$ with the addition of the possible new generators of length $n$, that is,
\begin{align*}
(N)_n=\{ x+y \,|\, x,y \in N, \ell(x+y)=n\} \cup (G)_n\: .
\end{align*}
Notice that, since $\ell(x+y)=\ell(x)+\ell(y)$,
\begin{align*}
\{ x+y \,|\, x,y \in N, \ell(x+y)=n\} &=\{ x+y \,|\, x\in (N)_i, y\in (N)_{j}, i+j=n\}\: ,
\end{align*}
and then,
\begin{align*}
|\{ x+y \,|\, x,y \in N, \ell(x+y)=n\}|=\sum_{\underset{0<i,j<n}{i+j=n}} |N|_i\cdot |N|_j\: .
\end{align*}
Since $G$ is a minimal generating set,
\begin{align*}
\{ x+y \,|\, x,y \in N, \ell(x)+\ell(y)=n\} \cap (G)_n=\emptyset\: .
\end{align*}
From there, the equality $|N|_n= |G|_n + \sum_{\underset{0<i,j<n}{i+j=n}} |N|_i\cdot |N|_j\: $.
\end{proof}

Let us show, through the above formula, that counting sequences of some finitely generated subgroupoids of $\M$ are related to \emph{Motzkin numbers}. A Motzkin path  of length $n$ is a path from $(0,0)$ to $(0,n)$ in the first octant of the planar grid $\mathbb{Z}^2$ with steps
\begin{tikzpicture}[scale=0.36]
\tikzstyle{vertex}=[circle, fill=black,minimum size=3pt,inner sep=1pt]
\node[vertex] (A_1) at (0,0) {};
\node[vertex] (A_2) at (1,1)   {};
\draw(A_1) --(A_2);
\end{tikzpicture}
,
\begin{tikzpicture}[scale=0.36]
\tikzstyle{vertex}=[circle, fill=black,minimum size=3pt,inner sep=1pt]
\node[vertex] (A_1) at (0,1) {};
\node[vertex] (A_2) at (1,0)   {};
\draw(A_1) --(A_2);
\end{tikzpicture}
, and
\begin{tikzpicture}[scale=0.36]
\tikzstyle{vertex}=[circle, fill=black,minimum size=3pt,inner sep=1pt]
\node[vertex] (A_1) at (0,0) {};
\node[vertex] (A_2) at (1,0)   {};
\draw(A_1) --(A_2);
\end{tikzpicture}
. Thus, Motzkin paths are a generalisation of Dyck paths, \cite{Donaghey1977}, \cite{Ainger1998}, \cite{Stanley2015}. Figure~\ref{MotzkinExamples} shows the nine possible Motzkin paths of length four.
The \emph{Motzkin number} $M_n$ counts the number of Motzkin paths of length $n$,
\cite[OEIS (A001006)] {Sloan2018TheEncy},
\begin{align*}
1, 1, 2, 4, 9, 21, 51, 127, 323, 835, 2188, 5798, 15511, 41835, 113634, 310572, \ldots
\end{align*}
for $n=0,1,2,\ldots$. Motzkin numbers satisfy the recurrence
$$M_{n+1}=M_{n}+\sum _{k=0}^{n-1}M_{k}M_{n-1-k}\: .$$
If we shift the sequence as $m_n=M_{n-1}$ we can define Motzkin numbers by the $\Cat$ transform in a compact form,
$$m_{n+1}=\Cat(m_n)\: .$$
These numbers are closely related to Catalan numbers through the binomial expressions
$$M_n=\sum_{k=0}^{\lfloor \frac{n}{2} \rfloor} {n \choose 2k} C_k\: , \qquad C_{n-1}=\sum_{k=0}^n {n \choose k} M_k\: .$$

\begin{example} \label{Example2unique} Consider the subgroupoid $2\M=\langle2 \rangle$. The minimal generating set has the sequence $|\{2\}|=0,1,0,0,0,0,\ldots$ and  the $\Cat$ transform yields
\begin{align*}
|2\M |_n=0,1, 0, 1, 0, 2, 0, 5, 0, 14, 0, 42, 0, 132, 0,1430, 0, 4862,  \ldots
\end{align*}
which is called the sequence of \emph{aerated Catalan numbers}. That sequence is related to the Motzkin numbers by the Euler transformation \cite{Donaghey1977}.
More in general, $|a\M|_n=c_{\frac{n}{\ell(a)}}$ if $\ell(a)|n$, and  $|a\M|_n=0$ otherwise. We already saw these sequences in Example~\ref{Exemple1}.
\end{example}

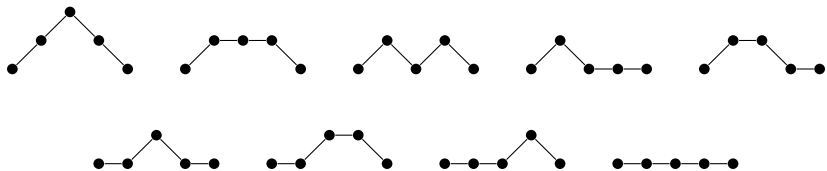
\begin{figure}[tb] 
\begin{tikzpicture}[scale=0.38]
\tikzstyle{vertex}=[circle, fill=black,minimum size=4pt,inner sep=1pt]
\node[vertex] (A_1) at (0,0) {};
\node[vertex] (A_2) at (1,1)   {};
\node[vertex] (A_3) at (2,2)  {};
\node[vertex] (A_4) at (3,1)  {};
\node[vertex] (A_5) at (4,0)  {};
\draw (A_1) -- (A_2) -- (A_3) --(A_4)--(A_5);

\node[vertex] (B_1) at (6,0) {};
\node[vertex] (B_2) at (7,1)   {};
\node[vertex] (B_3) at (8,1)  {};
\node[vertex] (B_4) at (9,1)  {};
\node[vertex] (B_5) at (10,0)  {};
\draw (B_1) -- (B_2) -- (B_3) --(B_4)--(B_5);

\node[vertex] (C_1) at (12,0) {};
\node[vertex] (C_2) at (13,1)   {};
\node[vertex] (C_3) at (14,0)  {};
\node[vertex] (C_4) at (15,1)  {};
\node[vertex] (C_5) at (16,0)  {};
\draw (C_1) -- (C_2) -- (C_3) --(C_4)--(C_5);

\node[vertex] (D_1) at (18,0) {};
\node[vertex] (D_2) at (19,1)   {};
\node[vertex] (D_3) at (20,0)  {};
\node[vertex] (D_4) at (21,0)  {};
\node[vertex] (D_5) at (22,0)  {};
\draw (D_1) -- (D_2) -- (D_3) --(D_4)--(D_5);

\node[vertex] (E_1) at (24,0) {};
\node[vertex] (E_2) at (25,1)   {};
\node[vertex] (E_3) at (26,1)  {};
\node[vertex] (E_4) at (27,0)  {};
\node[vertex] (E_5) at (28,0)  {};
\draw (E_1) -- (E_2) -- (E_3) --(E_4)--(E_5);
\end{tikzpicture}

\vspace{20pt}

\begin{tikzpicture}[scale=0.38]
\tikzstyle{vertex}=[circle, fill=black,minimum size=4pt,inner sep=1pt]
\node[vertex] (A_1) at (0,0) {};
\node[vertex] (A_2) at (1,0)   {};
\node[vertex] (A_3) at (2,1)  {};
\node[vertex] (A_4) at (3,0)  {};
\node[vertex] (A_5) at (4,0)  {};
\draw (A_1) -- (A_2) -- (A_3) --(A_4)--(A_5);

\node[vertex] (B_1) at (6,0) {};
\node[vertex] (B_2) at (7,0)   {};
\node[vertex] (B_3) at (8,1)  {};
\node[vertex] (B_4) at (9,1)  {};
\node[vertex] (B_5) at (10,0)  {};
\draw (B_1) -- (B_2) -- (B_3) --(B_4)--(B_5);

\node[vertex] (C_1) at (12,0) {};
\node[vertex] (C_2) at (13,0)   {};
\node[vertex] (C_3) at (14,0)  {};
\node[vertex] (C_4) at (15,1)  {};
\node[vertex] (C_5) at (16,0)  {};
\draw (C_1) -- (C_2) -- (C_3) --(C_4)--(C_5);

\node[vertex] (D_1) at (18,0) {};
\node[vertex] (D_2) at (19,0)   {};
\node[vertex] (D_3) at (20,0)  {};
\node[vertex] (D_4) at (21,0)  {};
\node[vertex] (D_5) at (22,0)  {};
\draw (D_1) -- (D_2) -- (D_3) --(D_4)--(D_5);
\end{tikzpicture}
\caption{The nine Motzkin paths of length four.}\label{MotzkinExamples}
\end{figure}

\begin{example} \label{Example23} By forbidding some figures in the path, we get a variant of Motzkin numbers \cite{Asinowski2020}. Consider the subgroupoid $\langle 2, 3_+\rangle$. The counting sequence of its generating set is $|\{2,3_+\}|_n=0,1,1,0,0,0,0,0,\ldots$. When we apply the $\Cat$ transform we obtain the sequence \cite[OEIS (A007477)] {Sloan2018TheEncy},
\begin{align*}
|\langle 2, 3_+\rangle |_n=0, 1, 1, 1, 2, 3, 6, 11, 22, 44, 90, 187, 392, 832, 1778, 3831, 8304,\ldots,
\end{align*}
which has generating function $(1-\sqrt{1-4x^2-4x^3})/2$, and it corresponds to the generating function of the number of Motzkin paths of length $n-2$ with the forbidden steps
\begin{tikzpicture}[scale=0.36]
\tikzstyle{vertex}=[circle, fill=black,minimum size=3pt,inner sep=1pt]
\node[vertex] (A_1) at (0,0) {};
\node[vertex] (A_2) at (1,0)   {};
\node[vertex] (A_3) at (2,1)  {};
\draw (A_1) --(A_2) -- (A_3);
\end{tikzpicture}
and
\begin{tikzpicture}[scale=0.35]
\tikzstyle{vertex}=[circle, fill=black,minimum size=3pt,inner sep=1pt]
\node[vertex] (A_1) at (0,0) {};
\node[vertex] (A_2) at (1,0)   {};
\node[vertex] (A_3) at (2,0)  {};
\draw (A_1) --(A_2)--(A_3);
\end{tikzpicture}
. See Figure~\ref{VariantMotzkinExamples1} for an example.
\end{example}

\begin{figure}[tb] 
\begin{tikzpicture}[scale=0.43]
\tikzstyle{vertex}=[circle, fill=black,minimum size=4pt,inner sep=1pt]
\node[vertex] (A_1) at (0,0) {};
\node[vertex] (A_2) at (1,1)   {};
\node[vertex] (A_3) at (2,2)  {};
\node[vertex] (A_4) at (3,1)  {};
\node[vertex] (A_5) at (4,0)  {};
\draw (A_1) -- (A_2) -- (A_3) --(A_4)--(A_5);

\node[vertex] (C_1) at (9,0) {};
\node[vertex] (C_2) at (10,1)   {};
\node[vertex] (C_3) at (11,0)  {};
\node[vertex] (C_4) at (12,1)  {};
\node[vertex] (C_5) at (13,0)  {};
\draw (C_1) -- (C_2) -- (C_3) --(C_4)--(C_5);

\node[vertex] (D_1) at (18,0) {};
\node[vertex] (D_2) at (19,1)   {};
\node[vertex] (D_3) at (20,1)  {};
\node[vertex] (D_4) at (21,0)  {};
\node[vertex] (D_5) at (22,0)  {};
\draw (D_1) -- (D_2) --(D_3);
\draw (D_3) -- (D_4) --(D_5);
\end{tikzpicture}
\caption{The three Motzkin paths of length four according to the variant exposed in Example~\ref{Example23}.}\label{VariantMotzkinExamples1}
\end{figure}
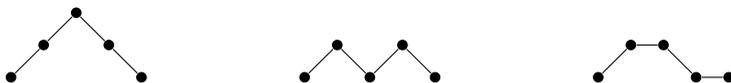

\begin{example} \label{Example233} Another variant of Motzkin paths permits decorating some steps with colours. We will draw the flat step with a double line,
\begin{tikzpicture}[scale=0.36]
\tikzstyle{vertex}=[circle, fill=black,minimum size=3pt,inner sep=1pt]
\node[vertex] (A_1) at (0,0) {};
\node[vertex] (A_2) at (1,0)   {};
\draw[thick,double](A_1) --(A_2);
\end{tikzpicture}
, to denote that it can be bicoloured independently of the other steps.
Consider the subgroupoid $\langle 2, 3_-, 3_+\rangle$. The counting sequence of its generating set is $|\{2,3_-,3_+\}|_n=0,1,2,0,0,0,0,0,\ldots$. By applying $\Cat$ we get the sequence \cite[OEIS (A253918)]{Sloan2018TheEncy},
\begin{align*}
|\langle 2, 3_-, 3_+\rangle |_n=0, 1, 2, 1, 4, 6, 12, 29, 56, 134, 300, 682, 1624, 3772, 9016, \ldots,
\end{align*}
which has generating function $(1-\sqrt{1-4x^2-8x^3})/2$, and it corresponds to the generating function of the number of Motzkin paths of length $n-2$ with the forbidden steps
\begin{tikzpicture}[scale=0.36]
\tikzstyle{vertex}=[circle, fill=black,minimum size=3pt,inner sep=1pt]
\node[vertex] (A_1) at (0,0) {};
\node[vertex] (A_2) at (1,0)   {};
\node[vertex] (A_3) at (2,1)  {};
\draw[thick,double] (A_1) --(A_2);
\draw (A_2) -- (A_3);
\end{tikzpicture}
and
\begin{tikzpicture}[scale=0.35]
\tikzstyle{vertex}=[circle, fill=black,minimum size=3pt,inner sep=1pt]
\node[vertex] (A_1) at (0,0) {};
\node[vertex] (A_2) at (1,0)   {};
\node[vertex] (A_3) at (2,0)  {};
\draw[thick,double] (A_1) --(A_2)--(A_3);
\end{tikzpicture}
. See Figure~\ref{VariantMotzkinExamples2} for an example of this variant of Motzkin paths.
\end{example}

\begin{figure}[tb] 
\begin{tikzpicture}[scale=0.43]
\tikzstyle{vertex}=[circle, fill=black,minimum size=4pt,inner sep=1pt]
\node[vertex] (A_1) at (0,0) {};
\node[vertex] (A_2) at (1,1)   {};
\node[vertex] (A_3) at (2,2)  {};
\node[vertex] (A_4) at (3,1)  {};
\node[vertex] (A_5) at (4,0)  {};
\draw (A_1) -- (A_2) -- (A_3) --(A_4)--(A_5);
\node[align=left] at (-1,1) {1 $\times$};

\node[vertex] (C_1) at (9,0) {};
\node[vertex] (C_2) at (10,1)   {};
\node[vertex] (C_3) at (11,0)  {};
\node[vertex] (C_4) at (12,1)  {};
\node[vertex] (C_5) at (13,0)  {};
\draw (C_1) -- (C_2) -- (C_3) --(C_4)--(C_5);
\node[align=left] at (8,1) {1 $\times$};

\node[vertex] (D_1) at (18,0) {};
\node[vertex] (D_2) at (19,1)   {};
\node[vertex] (D_3) at (20,1)  {};
\node[vertex] (D_4) at (21,0)  {};
\node[vertex] (D_5) at (22,0)  {};
\draw (D_1) -- (D_2);
\draw[thick,double] (D_2) --(D_3);
\draw (D_3) -- (D_4);
\draw[thick,double] (D_4) --(D_5);
\node[align=left] at (17,1) {4 $\times$};
\end{tikzpicture}
\caption{The six Motzkin paths of length four according to the variant exposed in Example~\ref{Example233}. Double lines can independently be bicoloured.}\label{VariantMotzkinExamples2}
\end{figure}
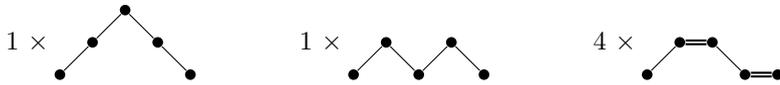

\begin{example} \label{Example33} Consider the subgroupoid $\langle 3_-, 3_+\rangle$. The counting sequence of its generating set is $|\{3_-,3_+\}|_n=0,0,2,0,0,0,0,0,\ldots$, and therefore,
\begin{align*}
|\langle 3_-, 3_+\rangle |_n=0, 0, 2, 0, 0, 4, 0, 0, 16, 0, 0, 80, 0, 0, 448, 0, 0, 2688, 0, 0, 16896,  \ldots
\end{align*}
with generating function $(1-\sqrt{1-8x^3})/2$. The sequence can be explicitly expressed as
\begin{align*}
|\langle 3_-, 3_+\rangle |_n=\begin{cases} 2^{\frac{n}{3}}c_{\frac{n}{3}} & \mbox{ if } 3|n, \\ 0 & \mbox{otherwise}.
\end{cases}
\end{align*}
\end{example}

Examples~\ref{Example2unique}, \ref{Example23}, and \ref{Example233} show that counting elements in a subgroupoid is closely related to Motzkin paths with multicoloured levels and avoiding some step patterns. However, the evidence is based on each case's coincidences of generating functions. That sets the problem of there being a general relation.

Return to the question of density, notice that all the examples viewed showed 0 or 1 densities, except for longitudinal subgroupoids, for which density is not well-defined.
In the style of the above examples, some numerical experiments make us think that finitely generated subgroupoids have always null density. The density tends to zero even for minimal generating sets with rapid counting sequences, such as $2^n$.
To finish, we would like to discover some intermediate-density subgroupoids.
To find such a subgroupoid, we need a minimal generating set that grows really fast. Consider a last example.

\begin{example} Consider the set $\M+1$. If we decompose any element $x \in \M+1$, we get $x=y+1$, for some $y \in \M$. However, $1\not \in \M+1$, meaning that $\M+1$ is a minimal generating set. Its counting sequence is the shifted Catalan sequence:
$$|\M+1|_n=0,c_1,c_2,c_3,c_4,c_5,c_6,  \ldots,$$
which has the generating function
$x(1-\sqrt{1-4x})/2$.
The $\Cat$ transform yields the sequence
\begin{align*}
|\langle \M+1 \rangle|_n=0, 1, 1, 3, 7, 21, 62, 197, 637, 2123, 7196, 24807, 86608, 305792, \ldots,
\end{align*}
which has generating function
\begin{align*}
\frac{1-\sqrt{1-2x+2x\sqrt{1-4x}}}{2}\: .
\end{align*}
We can approximately calculate the density as
\begin{align*}
\delta(\langle \M+1\rangle: \M) \approx 0.35361.
\end{align*}
We do not know the exact constant nor whether it is a rational number; we have calculated the above value evaluating the quotient $x_n=\frac{\Cat(|\M+1|_n)}{c_n}$ accelerated with the \emph{Aitken's $\Delta^2$ process} \cite{Shirley2017}, that is,
$\frac {x_n x_{n+2}-x_{n+1}^2}  {x_n+x_{n+2}-2 x_{n+1}}$
, for $n=5000$. However, the convergence velocity is so slow that, assuming it does converge, we can only fix a few digits.
In a very similar way, we can consider the minimal generating sets $\M+2$, or $\M+3_+$, for which we find the following approximated densities:
\begin{align*}
\delta(\langle \M+2\rangle: \M) & \approx 0.06683, \\
\delta(\langle \M+3_+\rangle: \M) & \approx 0.01588.
\end{align*}
\end{example}


%

\section*{Declarations}

\begin{itemize}

\item[] \textbf{Ethical Approval}. Not applicable. 
\item[] \textbf{Competing interests}. No conflict of interest exists.
\item[] \textbf{Authors' contributions}. Not applicable. 
\item[] \textbf{Funding}. Not applicable. 
\item[] \textbf{Availability of data and materials}. Not applicable.

\end{itemize}



\end{document}